\documentclass[12pt]{article}

\usepackage{amsfonts}
\usepackage{pslatex}
\usepackage{latexsym}
\newcommand{\proof}{\noindent{\it Proof: }}
\newcommand{\proofbox}{\mbox{ $\Box$}\\}

\newcommand{\R}{\mathbb{R}}

\begin{document}

\author{K\'aroly J. B\"or\"oczky\footnote{Supported by
OTKA grants 068398 and 049301, and by the EU Marie Curie TOK projects
DiscConvGeo and BudAlgGeo.} }

\title{Stability of the Blaschke-Santal\'o
and the affine isoperimetric inequality}

\newtheorem{lemma}{Lemma}[section]

\newtheorem{theo}[lemma]{Theorem}

\newtheorem{defi}[lemma]{Definition}

\newtheorem{coro}[lemma]{Corollary}

\newtheorem{conj}[lemma]{Conjecture}

\newtheorem{prop}[lemma]{Proposition}

\newtheorem{remark}[lemma]{Remark}

\newtheorem{example}[lemma]{Example}

\maketitle

\begin{center}
  Dedicated to Endre Makai on the occasion of his sixtieth birthday
\end{center}

\begin{abstract}
A stability version
of the Blaschke-Santal\'o inequality and
the affine isoperimetric inequality
for convex bodies of dimension $n\geq 3$
is proved. The first step is the reduction
to the case when the convex body is $o$-symmetric and
has axial rotational symmetry. This step works
for related inequalities compatible
with Steiner symmetrization.
Secondly, for these convex bodies, a stability version of the
characterization of ellipsoids by the fact that
each hyperplane section is centrally symmetric is established.
\end{abstract}


\section{Introduction}

Stability versions of geometric inequalities have been
investigated since the days of H. Minkowski, see the beautiful
survey of H. Groemer \cite{Gro93}, or K.J. B\"or\"oczky \cite{Bor05}
for some more recent results. Here we prove stability versions of
two classical inequalities originating from the beginning of the
20th century, the Blaschke-Santal\'o inequality and the affine
isoperimetric inequality. For all the basic affine invariant
notions, consult the thorough monograph of K. Leichtwei{\ss}
\cite{Lei98}, and for notions of convexity in general, see P.M.
Gruber \cite{Gru07} and R. Schneider \cite{Sch93}.

We write $o$ to denote the origin of $\R^n$,
$\langle \cdot,\cdot\rangle$ to denote
the standard scalar product, and $V(\cdot)$ to denote
volume. Let $B^n$ be the unit Euclidean ball
with volume $\kappa_n=V(B^n)$, and let $S^{n-1}=\partial B^n$.
A convex body $K$ in $\R^n$ is a compact convex set
with non--empty interior.
If $z\in{\rm int}K$,
then the polar of $K$ with respect to $z$ is the convex body
$$
K^z=\{x\in\R^n:\,\langle x-z,y-z\rangle\leq 1\mbox{ for any $y\in K$}\}.
$$
It is easy to see that $(K^z)^z=K$.
According to L.A. Santal\'o \cite{San49}
 (see also M. Meyer and A. Pajor \cite{MeP90}), there exists a unique
$z\in{\rm int}K$ minimizing the volume product
$V(K)V(K^z)$, which is called the Santal\'o point
of $K$. In this case $z$ is the centroid of $K^z$.
The celebrated  Blaschke-Santal\'o
inequality states that if $z$ is the  Santal\'o point
(or centroid) of $K$, then
\begin{equation}
\label{BS}
V(K)V(K^z)\leq\kappa_n^2,
\end{equation}
with equality if and only if $K$ is an ellipsoid.
The inequality was proved by W. Blaschke \cite{Bla17} for $n\leq 3$,
 and by L.A. Santal\'o \cite{San49} for all $n$.
The case of equality was characterized by
J. Saint-Raymond \cite{Sai81} among $o$-symmetric convex bodies,
and by C.M. Petty \cite{Pet85} among all convex bodies
(see also M. Meyer and A. Pajor \cite{MeP90},
D. Hug \cite{Hug96}, and M. Meyer and S. Reisner \cite{MeR06}
 for simpler proofs).

Our main task is to provide a stability version
of this inequality.
A natural tool is the  Banach-Mazur distance $\delta_{\rm BM}(K,M)$
of the convex bodies $K$ and $M$, which is defined by
$$
\delta_{\rm BM}(K,M)=\min\{\lambda\geq 1:\,
K-x\subset \Phi(M-y)\subset \lambda(K-x)\mbox{ for }
\Phi\in{\rm GL}(n),x,y\in\R^n\}.
$$
In particular, if $K$ and $M$ are $o$-symmetric, then
$x=y=o$ can be assumed. It follows
from a theorem of F. John \cite{Joh37}
that  $\delta_{\rm BM}(K,B^n)\leq n$
for any convex body $K$ in $\R^n$
(see also K.M. Ball \cite{Bal97}).

\begin{theo}
\label{BSstab}
If $K$ is a  convex body in $\R^n$, $n\geq 3$, with Santal\'o point $z$,
and
$$
V(K)V(K^z)>(1-\varepsilon)\kappa_n^2\mbox{ \ for $\varepsilon\in(0,\frac12)$,}
$$
then for some $\gamma>0$ depending only on $n$, we have
$$
\delta_{\rm BM}(K,B^n)<1+\gamma\,\varepsilon^{\frac1{6n}}
|\log \varepsilon|^{\frac16}.
$$
\end{theo}

Taking $K$ to be the convex body resulting from $B^n$
by cutting off two opposite caps of volume $\varepsilon$  shows that
 the exponent $1/(6n)$ cannot be replaced by anything larger
than $2/(n+1)$. Therefore the exponent of $\varepsilon$
is of the correct order. Since
$1/(6n)$ is most probably not the optimal
exponent of $\varepsilon$ in Theorem~\ref{BSstab}, no attempt was made to
find an explicit $\gamma$ in Theorem~\ref{BSstab}.
In principle, this would be possible following
the argument in this paper if
the exponent $1/(6n)$ is replaced by $1/(6n+6)$ (see the discussion
after (\ref{Brunnstab})).
We note that a stability version of the Blaschke-Santal\'o
inequality in the planar case is proved in
K.J. B\"or\"oczky, E. Makai,  M. Meyer, S. Reisner \cite{BMMR}, 
using a quite different method.

The literature about the Blaschke-Santal\'o inequality is
so extensive that only just a small portion can be discussed here.
The PhD thesis of K.M. Ball \cite{Bal} started off the quest for
suitable functional versions. This point of
view is for example pursued in M. Fradelizi and M. Meyer \cite{FrM07}
and S. Artstein, B. Klartag and V. D. Milman \cite{AKM04}.
Stability questions on a related problem are
discussed in M. Meyer and E. Werner \cite{MeW98}.

We note that the minimum of the volume product $V(K)V(K^z)$
is not known for convex bodies $K$ in $\R^n$ and $z\in K$
for $n\geq 3$. 
According to the well-known conjecture of
K. Mahler \cite{Mah39}, 
the volume product is minimized by simplices,
and among $o$-symmetric convex bodies by cubes.
The planar case was actually proved in \cite{Mah39}, and
simpler arguments are provided by M. Meyer \cite{Mey91} and 
M. Meyer and S. Reisner \cite{MeR06}.
For $n\geq 3$, the Mahler conjecture for $o$-symmetric convex bodies
 has been been verified
among unconditional bodies by J. Saint-Raymond \cite{Sai81}
(see also S. Reisner \cite{Rei87}),
and among zonoids by S. Reisner \cite{Rei86}
(see also Y. Gordon, M. Meyer and S. Reisner \cite{GMS88}).
The best lower bound for the volume product
of an $o$-symmetric convex body $K$ in $\R^n$ is
\begin{equation}
\label{MahlerKup}
V(K)V(K^o)>2^{-n}\kappa_n^2,
\end{equation}
due to G. Kuperberg \cite{GK}.
With a non-explicit
 exponential factor instead of $2^{-n}$,
it was proved by J. Bourgain and V.D. Milman \cite{BoM87}.

The Mahler conjecture for general convex bodies
was verified by
M. Meyer and S. Reisner \cite{MeR06} 
among polytopes of at most $n+3$ vertices.
In a yet unpublished revision
of \cite{GK}, G. Kuperberg also showed, based on (\ref{MahlerKup})
and the Rogers-Shephard inequality \cite{RoS58},
that
if $z\in{\rm int}K$ for a convex body $K$ in $\R^n$, then
\begin{equation}
\label{Mahlergen}
V(K)V(K^z)> 4^{-n}\kappa_n^2.
\end{equation}

It was probably W. Blaschke who first noticed that the
Blaschke-Santal\'o inequality is equivalent to the
affine isoperimetric inequality. This and other equivalent
formulations are discussed in depth in
E. Lutwak \cite{Lut93} and K. Leichtwei{\ss} \cite{Lei98}, Section 2.
To define the affine surface area of
a convex body $K$ in $\R^n$ , we always consider its
boundary endowed with the $(n-1)$-dimensional Hausdorff measure.
According to Alexandrov's theorem
(see P.M. Gruber \cite{Gru07}, p. 74),
$\partial K$ is twice differentiable in a generalized sense at almost
every point, hence the generalized Gau{\ss}-Kronecker curvature
$\kappa(x)$ can be defined at these $x\in\partial K$
(see K. Leichtwei{\ss} \cite{Lei98}, Section 1.2).
The affine surface area of $K$ is defined by
$$
\Omega(K)=\int_{\partial K}\kappa(x)^{\frac1{n+1}}dx.
$$
If $\partial K$ is $C^2$, then this definition is due to
W. Blaschke \cite{Bla16}.
Since then various equivalent definitions were given
for general convex bodies (including
the above one) by
K. Leichtwei{\ss} \cite{Lei86},
C. Sch\"utt and E. Werner \cite{ScW90} and
E. Lutwak \cite{Lut91}, 
which were shown to
be equivalent by C. Sch\"utt \cite{Sch92},
and G. Dolzmann and D. Hug \cite{DoH95}
(see K. Leichtwei{\ss} \cite{Lei98}, Section 2). The affine
surface area is a valuation invariant under volume preserving
affine transformations, and it is upper semi-continuous. These properties
are characteristic, as any upper semi-continuous
valuation on the space of convex bodies which is
invariant under volume preserving
affine transformations is a linear
combination of affine surface area, volume, and the Euler characteristic
by M. Ludwig and M. Reitzner \cite {LuM99}. We note that
affine surface area comes up e.g. in polytopal approximation
(see P.M. Gruber \cite{Gru07}, Section 11.2), 
in limit shape of lattice polygons
(see I. B\'ar\'any \cite{Bar97}), and many other applications
(see K. Leichtwei{\ss} \cite{Lei98}, Section 2).

The affine isoperimetric inequality states that
\begin{equation}
\label{affine}
\Omega(K)^{n+1}\leq \kappa_n^2 n^{n+1}V(K)^{n-1},
\end{equation}
with equality if and only if $K$ is an ellipsoid.
The inequality itself is due to  W. Blaschke  \cite{Bla16},
whose proof in $\R^3$ for convex bodies with $C^2$ boundaries
readily extends to general dimension and to 
general convex bodies.
W. Blaschke characterized
the equality case among convex bodies with $C^2$ boundary,
and this characterization was extended to all convex bodies
 by C.M. Petty \cite{Pet85}.
We note that W. Blaschke and L.A. Santal\'o
deduced the Blaschke-Santal\'o inequality from
the affine isoperimetric inequality. Here
we take a reverse path.

An inequality on p. 59 of E. Lutwak \cite{Lut91}
(see also Lemma 3.7 in D. Hug \cite{Hug96}, or (1106) in K. Leichtwei{\ss}
\cite{Lei98}) says that if $z\in{\rm int}K$, then
$$
\Omega(K)^{n+1}\leq  n^{n+1}V(K)^nV(K^z).
$$
Therefore Theorem~\ref{BSstab} yields

\begin{theo}
\label{affinestab}
If $K$ is a  convex body in $\R^n$, $n\geq 3$, and
$$
\Omega(K)^{n+1}>(1-\varepsilon)\kappa_n^2
n^{n+1}V(K)^{n-1}\mbox{ \ for $\varepsilon\in(0,\frac12)$},
$$
then for some $\gamma>0$ depending only on $n$, we have
$$
\delta_{\rm BM}(K,B^n)<1+\gamma\,\varepsilon^{\frac1{6n}}
|\log \varepsilon|^{\frac16}.
$$
\end{theo}

For convex bodies $K$ and $M$, we write $V_1(K,M)$ to denote the
mixed volume
$$
V_1(K,M)=\lim_{t\to 0}\frac{V(K+tM)-V(K)}{n\cdot t}
$$
(see T. Bonnesen and W. Fenchel \cite{BoF87}, Section 29,
or P.M. Gruber \cite{Gru07}, Section 6).
It satisfies $V_1(K,K)=V(K)$.
We write ${\cal K}^n_o$ to denote the family
of convex bodies whose centroid is $o$.
C.M. Petty \cite{Pet61} defined the geominimal
surface area by
$$
G(K)=\kappa_n^{-1/n}n\inf\{V_1(K,M^o)V(M)^{\frac1n}:M\in {\cal K}^n_o\}.
$$
It is also invariant under volume preserving affine transformations.
Positioning $K$ in a way such that $o$ is the Santal\'o point
of $K$ and taking $M=K^o$,  yields the so called
geominimal surface area inequality of C.M. Petty \cite{Pet74}
\begin{equation}
\label{geo}
G(K)\leq\kappa_n^{1/n}nV(K)^{\frac{n-1}n}, 
\end{equation}
 with equality if and only if $K$ is an ellipsoid.
From Theorem~\ref{BSstab} we directly obtain

\begin{theo}
\label{geostab}
If $K$ is a  convex body in $\R^n$, $n\geq 3$, and
$$
G(K)>(1-\varepsilon)
\kappa_n^{1/n}nV(K)^{\frac{n-1}n}\mbox{ \ for $\varepsilon\in(0,\frac12)$},
$$
then for some $\gamma>0$ depending only on $n$, we have
$$
\delta_{\rm BM}(K,B^n)<1+\gamma\,\varepsilon^{\frac1{6n}}
|\log \varepsilon|^{\frac16}.
$$
\end{theo}

One of our main tools is to reduce the proof of
Theorem~\ref{BSstab} to $o$-symmetric convex
bodies with axial rotational symmetry.

\begin{theo}
\label{rounding}
For any convex body $K$ in $\R^n$, $n\geq 2$, with
$\delta_{\rm BM}(K,B^n)>1+\varepsilon$ for $\varepsilon>0$,
there exists an $o$-symmetric convex body $C$ with
axial rotational symmetry and a constant
$\gamma>0$ depending on $n$
such that 
$\delta_{\rm BM}(C,B^n)>1+\gamma\varepsilon^2$ and
$C$ results from $K$ as a limit
of subsequent Steiner symmetrizations and
affine transformations.
\end{theo}
{\bf Remark: } If $K$ is $o$-symmetric, then
$1+\gamma\varepsilon^2$ can be replaced
by $1+\gamma\varepsilon$.  In particular, if $K$ is
$o$-symmetric, then wherever the factor $1/6$ occurs
in Theorems~\ref{BSstab}, \ref{affinestab} and \ref{geostab},
it can be replaced by $1/3$. \\

Theorem~\ref{rounding} shows that it is possible
to use Steiner-symmetrization to obtain a convex body that
is highly symmetric but still far from being an ellipsoid.
On the other hand, B. Klartag \cite{Kla04} proved
that any convex body $K$ in $\R^n$ gets $\varepsilon$ close
to some ball after suitable chosen $cn^4|\log \varepsilon|^2$
Steiner symmetrizations where $c>0$ is an absolute constant.

After discussing the basic tools 
such as the isotropic
position of convex bodies
and Steiner symmetrization in Section~\ref{sectools},
we prove Theorem~\ref{rounding} in Section~\ref{secrounding}.
A stability version of the False Centre theorem
in a special case is presented in Section~\ref{secFalsestab}, which
combined with  Theorem~\ref{rounding} leads
to the proof of Theorem~\ref{BSstab} in Section~\ref{secBSstab}.
For stability versions
of some other classical geometric characterizations
of ellipsoids, see, e.g., H. Groemer \cite{Gro94}
and P.M. Gruber \cite{Gru97}.

\section{Some tools}
\label{sectools}

\subsection{Isotropic position}

In this paper, we use the term isotropic position
in a weak sense.
More precisely, we say that a convex body $K$ in $\R^n$ is
in weak isotropic position if
its centroid is $o$, and $\int_K\langle u,x\rangle^2dx$ is
independent of $u\in S^{n-1}$.
In particular, in this case
$$
\int_K\langle u,x\rangle^2dx=L_K^2 V(K)^{\frac{n+2}n}
$$
for any $u\in S^{n-1}$
(see, e.g.,  A.A. Giannopoulos and V.D. Milman \cite{GiM04}), and
the Legendre ellipsoid (the ellipsoid of inertia)  is
a ball.
For any convex body $C$ there is a volume
preserving affine transformation $T$ such that
$TC$ is in weak isotropic position. In the
literature, two diffferent normalizations
are used. Either $V(K)=1$ (see, e.g., 
A.A. Giannopoulos and V.D. Milman \cite{GiM04}), or
$\|v\|^2=\int_K\langle v,x\rangle^2dx$ for
any $v\in\R^n$ (see, e.g., 
R. Kannan, L. Lov\'asz and M. Simonovits \cite{KLS95}).
In this paper, if $K$ is in weak isotropic position,
then we compare it to balls,
therefore we frequently assume $V(K)=\kappa_n$.

It is known that $L_K$ is minimized by ellipsoids (see F. John \cite{Joh37}
or A.A. Giannopoulos and V.D. Milman \cite{GiM04}).
It follows by Gy. Sonnevend \cite{Son89} (see also
R. Kannan, L. Lov\'asz and M. Simonovits \cite{KLS95})
that if $K$ is in weak isotropic position, then
$$
K\subset L_K^{\frac2{n+2}}V(K)^{\frac1n}\sqrt{n(n+2)} B^n.
$$
Now $L_K\leq c_0\sqrt[4]{n}$ for some absolute constant $c_0$ according
to  B. Klartag \cite{Kla06}. Therefore,
if $V(K)=\kappa_n$ and $K$ is in weak
isotropic position, then
\begin{equation}
\label{isodiam}
K\subset c\sqrt{n}\,B^n
\end{equation}
for some absolute constant $c\geq 1$.

For properties of $o$-symmetric convex bodies in isotropic
position,
see the discussion in V.D. Milman and A. Pajor \cite{MiP89}.

\subsection{Steiner symmetrization}
\label{secSteiner}

Given a convex body $K$ in $\R^n$ and a hyperplane
$H$, for any $l$ orthogonal to $H$ and intersecting
$K$, translate $l\cap K$ along $l$ in a way such that
the midpoint of the image lies in $H$. The union
of these images is the Steiner symmetrial
$K_H$ of $K$ with respect to $H$.
It follows that $K_H$ is convex, $V(K_H)=V(K)$,
and, if the centroid of $K$ lies in $H$, then
it coincides with the centroid of $K_H$.

We write $|\cdot|$ to denote the $(n-1)$-dimensional Lebesgue measure,
where the measure of the empty set is defined to be zero.
For $u\in S^{n-1}$  and $t\in\R$, let
$u^\bot$ denote the linear $(n-1)$-space orthogonal to $u$, let
$h_K(u)=\max_{x\in K}\langle u,x\rangle$ be the support function
of $K$, and let
$$
K(u,t)=K\cap (tu+u^\bot).
$$
If $M$ is a compact convex set of dimension $n-1$,
then the classical Brunn-Minkowski inequality
(see, e.g., T. Bonnesen and W. Fenchel \cite{BoF87}, p. 94,
 P.M. Gruber \cite{Gru07},
Section 8.1,
or the monograph R. Schneider \cite{Sch93}, which is
solely dedicated to the Brunn-Minkowski theory) yields
\begin{equation}
\label{Brunn}
|\mbox{$\frac12$}(M-M)|\geq |M|.
\end{equation}

K.M. Ball proved in his PhD
thesis \cite{Bal} that Steiner symmetrization
through $u^\bot$ for $u\in S^{n-1}$
increases $V(K^o)$ if $K$ is $o$-symmetric.
The basis of his argument is the observation that
for $\widetilde{K}=K_{u^\bot}$, we have
\begin{equation}
\label{section}
\mbox{$\frac12$}\left(K^o(u,t)-K^o(u,t)\right)\subset \widetilde{K}^o(u,t)-tu
\end{equation}
(see also M. Meyer and A. Pajor \cite{MeP90}).
Here the $(n-1)$-measure of the
 left hand side is at least $|K^o(u,t)|$ according
to the Brunn-Minkowski inequality, hence the Fubini Theorem
yields $V(\widetilde{K}^o)\geq V(K^o)$.
K.M. Ball's result was further exploited by
 M. Meyer and A. Pajor \cite{MeP90}. 
The ideas and statements in \cite{MeP90} yield the following.

\begin{lemma}[Meyer,Pajor]
\label{BSsteiner}
Let $K$ be a convex body in $\R^n$, and let $H$ be
a hyperplane. If $z$ and $z'$ denote the Santal\'o
points of $K$ and $K_H$ , respectively, then
$z'\in H$, and $V(K^z)\leq V((K_H)^{z'})$.
\end{lemma}

This statement is more explicit in Theorem~1 of
M. Meyer and S. Reisner \cite{MeR06}
(see the proof of Theorem~13 in \cite{MeR06}).

\section{Proof of Theorem~\ref{rounding}}
\label{secrounding}

The following lemma is the basis of the
proof of Theorem~\ref{rounding}.

\begin{lemma}
\label{rounding0}
Let $K$ be a convex body in $\R^n$.
If $\delta_{\rm BM}(K,B^n)>1+\varepsilon$ for $\varepsilon>0$,
then there exists a convex body $C$ with
axial rotational symmetry that
results from $K$ as a limit of
 subsequent Steiner symmetrizations and
affine transformations, and satisfies
$\delta_{\rm BM}(C,B^n)>1+\gamma\varepsilon$,
where $\gamma>0$ depends only on $n$.
Moreover if $K$ is $o$-symmetric, then so is $C$.
\end{lemma}
\proof We may assume that $V(K)=\kappa_n$ and $K$
is in weak isotropic position.
Using the $c\geq 1$ from (\ref{isodiam}), we claim that
 there exists some
$u\in S^{n-1}$ such that
\begin{description}
\item{(i) } either
$h_K(u)\geq 1+\frac{\varepsilon}4$ and
$V(K\backslash B^n)\leq\tilde{\gamma}\varepsilon$
for $\tilde{\gamma}=\frac1{4c^2n}\int_{B^n}\langle u,x\rangle^2dx$,
\item{(ii) } or $h_K(u)\leq 1-\frac{\tilde{\gamma}}{n\kappa_n}\,\varepsilon$.
\end{description}
To prove this statement, let $h_K$ attain
its maximum on $S^{n-1}$ at $v\in S^{n-1}$, and its minimum at $w\in S^{n-1}$.
If $h_K(w)\leq 1-\frac{\varepsilon}4$, then $u=w$ works, thus we 
may assume $h_K(w)\geq 1-\frac{\varepsilon}4$.
Since $\delta_{\rm BM}(K,B^n)>1+\varepsilon$,
it follows that $h_K(v)\geq 1+\frac{\varepsilon}4$.
Now if $V(K\backslash B^n)\leq\tilde{\gamma}\varepsilon$,
then we are done again, hence we may assume
$V(B^n\backslash K)=V(K\backslash B^n)\geq\tilde{\gamma}\varepsilon$.
We conclude that
$h_K(w)\leq 1-\frac{\tilde{\gamma}}{n\kappa_n}\,\varepsilon$,
which completes
the proof of (i) and (ii).

Let $C$ be the
image of $K$
after applying first Schwarz rounding
(see P.M. Gruber \cite{Gru07}, p. 178)
in the direction of $u$, and secondly
the linear transformation that
 dilates by the factor $h_K(u)^{-1}$ in the direction
of $u$, and by the factor $h_K(u)^{\frac1{n-1}}$ orthogonal
to $u$. Since Schwarz rounding can be obtained
as the limit of repeated applications of Steiner symmetrizations
through hyperplanes containing the line $\R u$,
we have $V(C)=V(K)$ and
$o$ is the centroid of $C$ (see Section~\ref{secSteiner}). 
The linear transformation
following the Schwarz rounding ensures that $u\in \partial C$.

Let $h=h_K(u)$ and $\tilde{h}=h_K(-u)$. In the case of (ii),
$L_K\geq L_{B^n}$ yields
\begin{eqnarray}
\nonumber
\int_C\langle u,x\rangle^2dx&=&\int_0^1r^2|C(u,r)|\,dr+\int_0^{\tilde{h}/h}r^2|C(-u,r)|\,dr\\
\nonumber
&=& \int_0^1r^2h|K(u,hr)|\,dr+\int_0^{\tilde{h}/h}r^2h|K(-u,hr)|\,dr\\
\nonumber
&=&\frac1{h^2}\left(\int_0^hs^2|K(u,s)|\,ds+\int_0^{\tilde{h}}s^2|K(u,s)|\,ds\right)\\
\nonumber
&=&\frac1{h^2}\int_K\langle u,x\rangle^2dx\\
\nonumber
&=&\frac1{h^2}L_K^2 \kappa_n^{\frac{n+2}n}\\
\nonumber
&\geq & \frac1{h^2}\int_{B^n}\langle u,x\rangle^2dx\\
\label{(ii)cons}
&>& (1+\mbox{$\frac{\tilde{\gamma}}{n\kappa_n}$}\,\varepsilon)
\int_{B^n}\langle u,x\rangle^2dx.
\end{eqnarray}

In the case of (i), we have
$K\subset c\sqrt{n}B^n$ according to (\ref{isodiam}).
It follows that
\begin{eqnarray}
\nonumber
\int_C\langle u,x\rangle^2dx
&=&\frac1{h^2}\int_K\langle u,x\rangle^2dx\\
\nonumber
&<&\frac1{h^2}\left(\int_{B^n}\langle u,x\rangle^2dx+
c^2nV(K\backslash B^n)\right)\\
\nonumber
&\leq & \frac{1+\frac{\varepsilon}4}{h^2}
\int_{B^n}\langle u,x\rangle^2dx\\
\label{(i)cons}
&<& (1-\mbox{$\frac{\varepsilon}8$})\int_{B^n}\langle u,x\rangle^2dx.
\end{eqnarray}

Let $\delta_{\rm BM}(C,B^n)=1+\eta$, where we may assume that $\eta\in(0,1)$. Since $C$ has axial rotational symmetry 
around $\R u$, and $o$ is the centroid of $C$, there exist
$\gamma_1>0$
depending only on $n$, and an $o$-symmetric ellipsoid $E$ 
with axial rotational symmetry 
around $\R u$ such that $E\subset C\subset (1+\gamma_1\eta)E$.
It follows by $V(C)=V(B^n)$ and $u\in\partial C$ that 
there exists a 
$\gamma_2>0$ depending only on $n$ such that
$$
(1+\gamma_2\eta)^{-1}B^n\subset C\subset (1+\gamma_2\eta)B^n.
$$
Therefore, we conclude Lemma~\ref{rounding0}
by (\ref{(i)cons}) in the case of (i),
and by (\ref{(ii)cons})  in the case of (ii).
\proofbox

Let us write $W(M)$ to denote the mean width of a
planar compact convex set $M$. In particular $\pi W(M)$
is the perimeter of $M$. Writing $R(M)$ and
$r(M)$ to denote the circum- and the inradius of $M$,
and $A(M)$ to denote the area of $M$,
the Bonnesen inequality (appearing
first in W. Blaschke \cite{Bla22},
see H. Groemer \cite{Gro93} for more references) states
\begin{equation}
\label{Bonnesen}
W(M)^2-\mbox{$\frac4{\pi}$}\,A(M)\geq (R(M)-r(M))^2.
\end{equation}

To prove Theorem~\ref{rounding} for convex bodies in $\R^n$,
we need the following statement.

\begin{prop}
\label{planarSteiner}
If $M$ is a planar compact convex set in $\R^2$
with an axis of symmetry satisfying
$\delta_{\rm BM}(K,B^2)>1+\varepsilon$ for $\varepsilon\in(0,1)$,
then there exist orthogonal lines $l_1$ and $l_2$ such that
$\delta_{\rm BM}((K_{l_1})_{l_2},B^2)>1+c'\varepsilon^2$
for $c'=0.001$.
\end{prop}
\proof Let $l$ be the line of symmetry of $K$.
We may assume that $A(K)=\pi$, and that $l$ intersects
$K$ in a segment of length $2$ whose midpoint is $o$.

First we try Steiner symmetrization through
$l$, and the line $l'$ that is orthogonal
to $l$ through $o$. 
If $\delta_{\rm BM}((K_{l})_{l'},B^2)> 1+c'\varepsilon^2$,
then we are done. Otherwise there is an ellipse $E$ whose
principal axis are
contained in $l$ and $l'$ such that
$$
E\subset (K_{l})_{l'}\subset (1+c'\varepsilon^2)E.
$$
We deduce that
\begin{equation}
\label{korkozel}
(1+c'\varepsilon^2)^{-3}B^2\subset (K_{l})_{l'}.
\end{equation}
Since $\delta_{\rm BM}(K,B^2)>1+\varepsilon$,
it follows that $R(K)-r(K)\geq \varepsilon/2$.
Therefore, the Bonnesen inequality (\ref{Bonnesen})
yields
$$
W(K)\geq 2\cdot\left(1+\frac{\varepsilon^2}{16}\right)^{\frac12}.
$$
In particular, if the distance of $x_1,x_2\in \partial K$
is the diameter of $K$, then
$$
\|x_1-x_2\|> 2\cdot(1+c'\varepsilon^2)^5.
$$
Next let $s$ be the  segment orthogonal to $x_1-x_2$
and of length $2(1+c'\varepsilon^2)^{-3}$. Since $K$ is symmetric
through $l$, (\ref{korkozel}) yields that $s'\subset K$
for a translate $s'$ of $s$. We deduce that
the convex hull $Q$ of $x_1,x_2$ and $s'$ satisfies
$$
A(Q)> 2\cdot (1+c'\varepsilon^2)^2.
$$
Let $l_1$ be the line determined by $x_1$ and $x_2$,
let $l_2$ be an orthogonal line, and let $K'=(K_{l_1})_{l_2}$.
Then $K'$ contains a quadrilateral of area larger than
$(1+c'\varepsilon^2)^2 \cdot\frac2{\pi}\,A(K')$,
which in turn yields that
$\delta_{\rm BM}(K',B^2)>1+c'\varepsilon^2$.
\proofbox

\noindent{\it Proof of Theorem~\ref{rounding}: }
If $K$ is $o$-symmetric, then Lemma~\ref{rounding0}
yields Theorem~\ref{rounding}.
Even if $K$ is not $o$-symmetric, we may assume that
$K$ has rotational symmetry around $\R u$
for some $u\in S^{n-1}$ and
satisfies $\delta_{\rm BM}(K,B^n)>1+\gamma\varepsilon$
for the $\gamma$ in Lemma~\ref{rounding0}. We deduce by
Proposition~\ref{planarSteiner} that there exist
orthogonal hyperplanes $H_1$ and $H_2$ containing $o$ such that
$H_1\cap H_2$ is orthogonal to $u$, and
$\delta_{\rm BM}(\widetilde{K},B^n)>1+c'\gamma^2\varepsilon^2$
for $\widetilde{K}=(K_{H_1})_{H_2}$
and the absolute constant $c'$ of Proposition~\ref{planarSteiner}.
Since $\widetilde{K}$ is $o$-symmetric,
the $o$-symmetric case of Lemma~\ref{rounding0}
applied to $\widetilde{K}$ yields Theorem~\ref{rounding}
for $K$.
\proofbox

\section{Stability of the False Centre Theorem
in a special case}
\label{secFalsestab}

For any convex body $K$ in $\R^n$,
P.W. Aitchison, C.M. Petty, C.A. Rogers \cite{APR71}
and D.G. Larman \cite{Lar74}
proved the False Centre Theorem, which states that
if there exists a point $p$ such that all hyperplane
sections of $K$ by hyperplanes passing through $p$ are centrally symmetric,
then $K$  is either symmetric through $p$ or an ellipsoid.
An important part of their proof is concerned
with the case when $K$ is $o$-symmetric and has
axial rotational symmetry.  We will deal with this
special case in Lemma~\ref{Falsestab}.

We measure how close a compact convex set $M$
is to be centrally symmetric by the so called
Minkowski measure of symmetry $q(M)$. It is defined by
(see, e.g., B. Gr\"unbaum \cite{Gru63})
$$
q(M)=\min\{\lambda\geq 1: \exists x\in M, -(M-x)\subset\lambda (M-x)\}.
$$
Obviously, $q(M)=1$ if and only if $M$ is centrally symmetric.
Moreover, it is known essentially since 
the time of H. Minkowski that $q(M)\leq
n$ for $M\subset \R^n$, where equality holds only for
$n$-dimensional simplices. To prove Lemma~\ref{Falsestab}, we need
the following estimate:

\begin{prop}
\label{Minksym}
Let $g$ be a positive concave function on $(-\varrho,\varrho)$ 
for $\varrho>0$,
and let $M$ be the compact convex set that is the closure
of the convex hull of the graphs of $g$ and $-g$. 
If $q(M)\leq 1+\varepsilon$
for $\varepsilon>0$, then for any $t\in(0,\varrho)$, we have
$$
\left(1+\frac{2\varrho\varepsilon}{\varrho-t}\right)^{-1} g(t)
\leq g(-t)\leq \left(1+\frac{2\varrho\varepsilon}{\varrho-t}\right) g(t).
$$
\end{prop}
\proof We may assume that $\varrho=1$.
Writing $u$ to denote the first coordinate
unit vector, the condition $q(M)\leq 1+\varepsilon$ yields that
$M\subset-(1+\varepsilon)M+pu$, where $|p|\leq \varepsilon$.
In particular, for any $t\in(0,1)$, we have
$$
g(-t)\leq (1+\varepsilon)g(\mbox{$\frac{t+p}{1+\varepsilon}$}).
$$
If $\frac{t+p}{1+\varepsilon}\geq t$, then,
considering the points $(-1,0)$, $(t,g(t))$ and 
$(\frac{t+p}{1+\varepsilon},g(\frac{t+p}{1+\varepsilon}))$
of $\partial M$, leads to
$$
g(\mbox{$\frac{t+p}{1+\varepsilon}$})\leq
\frac{1+\frac{t+p}{1+\varepsilon}}{1+t}\cdot g(t)\leq
\frac1{1+\varepsilon}\cdot\left(1+\frac{2\varepsilon}{1+t}\right)
\cdot g(t),
$$
and if $\frac{t+p}{1+\varepsilon}\leq t$, then
$$
g(\mbox{$\frac{t+p}{1+\varepsilon}$})\leq
\frac{1-\frac{t+p}{1+\varepsilon}}{1-t}\cdot g(t)\leq
\frac1{1+\varepsilon}\cdot\left(1+\frac{2\varepsilon}{1-t}\right)
\cdot g(t).
$$
In turn, we conclude the required upper bound
for $g(-t)$. To get the lower bound,
one applies the same argument for $h(t)=g(-t)$.
\proofbox

\begin{lemma}
\label{Falsestab}
Let $K$ be an $o$-symmetric convex body in $\R^n$, $n\geq 3$,
with axial rotational symmetry. If $\delta_{\rm BM}(K,B^n)>1+\varepsilon$
for $\varepsilon>0$,
then there exists a hyperplane $H$ intersecting $\frac23 K$ such that
 $q(H\cap K)\geq 1+c_1\varepsilon^3|\log\varepsilon|^{-1}$,
 where $c_1>0$ is an absolute constant.
\end{lemma}
{\rm Remark: } In the proof we only use hyperplanes that pass
 through one of the endpoints of the axis of $K$, therefore
we do have a stability version of the False Centre Theorem in this
special case. We believe that in Lemma~\ref{Falsestab}, the term $\varepsilon^3|\log\varepsilon|^{-1}$
can be improved to $\varepsilon$.\\
\proof  The proof is based on ideas of
P.W. Aitchison, C.M. Petty, C.A. Rogers \cite{APR71}.
We may assume that $u,-u\in\partial K$
where $u\in S^{n-1}$ and $\R u$ is the axis of symmetry of $K$.
We prove the lemma in the following form.
There exists a positive absolute constant $\tilde{c}$
such that if
$\varepsilon\in(0,4^{-4})$ and
 $q(H\cap K)\leq 1+\varepsilon$
holds for any hyperplane $H$ intersecting 
$\frac23 K$ and containing $-u$, then 
$\delta_{\rm BM}(K,B^n)\leq 
1+\tilde{c}\varepsilon^{\frac13}|\log \varepsilon|$.
To prove this statement, we may assume that
$\partial K$ is $C^1$.

By the symmetry of $K$, we may assume that $n=3$.
Let $v\in S^2$ be orthogonal to $u$,
and let $L$ be the linear plane spanned by $u$ and $v$.
There exists a non-negative even concave function $r$ on $[-1,1]$
such that $tu+r(t)v\in\partial K$ for $t\in[-1,1]$ and $r(0)=1$.
This $r$ is differentiable on $(-1,1)$
because $\partial K$ is $C^1$.
To prove that $K$ is close to some ellipsoid is equivalent
 to showing that the function
$$
f(t)=\frac{1-r(t)^2}{t^2}
$$
is essentially the constant one function on $(0,1)$.
In this proof, the implied constant in $O(\cdot)$
is always some absolute constant.

For $m\in(0,\frac14]$, let $H$ be the plane
containing $-u$ and $(1-m)u+r(1-m)v$, whose normal vectors
are contained in $L$, and let $\eta=\frac{r(1-m)}{2-m}$
be the ``slope'' of $H\cap L$.
In particular, if $l\subset H$ is a line
orthogonal to $L$ and passing through
the point $tu+\eta(1+t)\,v$, $t\in(-1,1-m)$, then
it intersects $K$ in a segment of length
$2\sqrt{r(t)^2-\eta^2(1+t)^2}$.
Since $q(H\cap K)\leq 1+\varepsilon$, Proposition~\ref{Minksym}
yields for any $t\in[0,1-m)$ that
$$
r(-t-m)^2-\eta^2(1-t-m)^2
\left\{
\begin{array}{cl}
\leq &
\left(1+\frac{(2-m)\varepsilon}{1-m-t}\right)^2
(r(t)^2-\eta^2(1+t)^2)\\
\geq &
\left(1+\frac{(2-m)\varepsilon}{1-m-t}\right)^{-2}
(r(t)^2-\eta^2(1+t)^2)
\end{array}\right.
$$
In particular, if $t\in[0,1-2m]$, then
\begin{equation}
\label{rtdiff}
r(t)^2-r(t+m)^2=\eta^2(2t+m)(2-m)+O\left(\frac{\varepsilon}{1-t}\right).
\end{equation}
For $t=0$, we have
\begin{equation}
\label{r0diff}
m^2 f(m)=1-r(m)^2=\eta^2m(2-m)+O(\varepsilon)=
\frac{m\cdot r(1-m)^2}{2-m}+O(\varepsilon).
\end{equation}
If $t\in[m,1-2m]$, then (\ref{rtdiff}) can be written in the form
\begin{eqnarray*}
(t+m)^2f(t+m)-t^2f(t)
&=&\eta^2(2t+m)(2-m)+O\left(\frac{\varepsilon}{1-t}\right)\\
&=&(2tm+m^2) f(m)+O\left(\frac{(2t+m)\varepsilon}{m}\right),
\end{eqnarray*}
therefore,
\begin{equation}
\label{ft+m}
 f(t+m)=\frac{t^2}{(t+m)^2}\,f(t)+\frac{2tm+m^2}{(t+m)^2}\,f(m)
+O\left(\frac{\varepsilon}{mt}\right).
\end{equation}
We deduce by (\ref{ft+m}) and induction
that if $i=2,\ldots,\lfloor\frac1m-1\rfloor$, then
\begin{equation}
\label{fim}
 f(im)=f(m)+O\left(\sum_{j=1}^{i-1}\frac{\varepsilon}{jm^2}\right)
=f(m)+O\left(\frac{\varepsilon|\log m|}{m^2}\right).
\end{equation}

We define
$$
\tilde{m}=\frac1{4\lfloor\varepsilon^{-\frac13}|
\log \varepsilon|^{-\frac13}\rfloor}.
$$
By definition, $\tilde{m}$ satisfies
\begin{equation}
\label{mtilde}
\frac{\varepsilon|\log \tilde{m}|}{\tilde{m}^2}=O(\tilde{m})
\mbox{ \ and \ }\tilde{m}\leq \frac18.
\end{equation}
We claim that
\begin{equation}
\label{ftildem}
f(\tilde{im})=1+O(\tilde{m}) \mbox{ \ for \ $i=1,\ldots,\frac1{\tilde{m}}-1$}.
\end{equation}
First we observe that according to (\ref{fim}), (\ref{mtilde}),
 and the definition
of $f$, we have
$$
f(i\tilde{m})=f(1-\tilde{m})+O(\tilde{m})\leq (1-\tilde{m})^{-2}+O(\tilde{m})=
1+O(\tilde{m})
$$
 for $i=1,\ldots,\frac1{\tilde{m}}-1$.
On the other hand, it follows by (\ref{r0diff}) that
$$
r(1-\tilde{m})^2=(2-\tilde{m})\tilde{m}\,f(\tilde{m})
+O\left(\frac{\varepsilon}{\tilde{m}}\right)=O(\tilde{m}).
$$
In particular,
$$
f(1-\tilde{m})=\frac{1-r(1-\tilde{m})^2}{(1-\tilde{m})^2}\geq
\frac{1-O(\tilde{m})}{(1-\tilde{m})^2}\geq 1-O(\tilde{m}),
$$
which in turn yields (\ref{ftildem}) by  (\ref{fim}) and (\ref{mtilde}).

Finally we verify that if $\tilde{m}\leq t\leq 1-\tilde{m}$, then
\begin{equation}
\label{ft}
f(t)=1+O(\tilde{m}) \mbox{ \ for $t\in[\tilde{m},1-\tilde{m}]$}.
\end{equation}
First let $t\in[\frac12,1-\tilde{m}]$. In this case,
$$
f'(t)=\frac{-2r(t)r'(t)t-2(1-r(t)^2)}{t^3}\geq -16
$$
as $r'(t)\leq 0$.
Since there exists an integer $i\leq \frac1{\tilde{m}}-2$
such that $\frac12\leq i\tilde{m}\leq t\leq (i+1)\tilde{m}$,
we deduce (\ref{ft}) from (\ref{fim}) and (\ref{mtilde}).

Next let $t\in[\tilde{m},\frac12]$. There exist integers $j$ and $i$ such that
$m\in [\tilde{m},2\tilde{m}]$ holds for $m=t/j$, and $im\in[\frac12,1-m]$,
thus (\ref{fim}) and the previous case of (\ref{ft}) yield
$$
 f(t)=f(jm)=f(m)+O\left(\frac{\varepsilon|\log m|}{m^2}\right)
=f(im)+O\left(\frac{\varepsilon|\log m|}{m^2}\right)
= 1+O(\tilde{m}).
$$
With this, we have proved (\ref{ft}), 
which in turn yields Lemma~\ref{Falsestab}.
\proofbox
 
From Lemma~\ref{Falsestab} we immediately obtain:

\begin{coro}
\label{Falsestabcor}
Let $K$ be an $o$-symmetric convex body in $\R^n$
with axial rotational symmetry. If $\delta_{\rm BM}(K,B^n)>1+\varepsilon$
for $\varepsilon>0$,
then there exist $u\in S^{n-1}$ and $a\in (0,\frac23)$
such that
$$
q(K(u,h_K(u)t))\geq 1+c_2\varepsilon^3|\log\varepsilon|^{-1}
\mbox{ \ for $t\in(a,a+c_2\varepsilon^3|\log\varepsilon|^{-1})$},
$$
where $c_2>0$ is an absolute constant.
\end{coro}

\section{Proof of Theorem~\ref{BSstab}}
\label{secBSstab}

We will need a stability version of the Brunn-Minkowski inequality.
According to V.I. Diskant \cite{Dis73}, if
$M$ is a compact convex set of dimension $n-1$ with
$q(M)\geq 1+\tau$, then
\begin{equation}
\label{Brunnstab}
|\mbox{$\frac12$}(M-M)|\geq (1+\gamma\tau^{n-1})|M|,
\end{equation}
for $\gamma>0$ depending on $n$ (see also H. Groemer \cite{Gro93}).
Here no explicit $\gamma$ is known. Actually H. Groemer \cite{Gro88}
proved a stability estimate with explicit $\gamma$ but with the exponent
$n$ instead of $n-1$.

In this section, $\gamma_1,\gamma_2,\ldots$ denote
positive constants depending only on $n$.
We prove Theorem~\ref{BSstab} in the following equivalent form:
If $K$ is a convex body in $\R^n$ with Santal\'o point $z$
and $\delta_{\rm BM}(K,B^n)>1+\varepsilon$ for $\varepsilon>0$,
then (\ref{BSstabeq}) holds.

It follows from Theorem~\ref{rounding} and Lemma~\ref{BSsteiner}
that there exists an $o$-symmetric convex body $C$ with
axial rotational symmetry such that
$\delta_{\rm BM}(C,B^n)>1+\gamma_1\varepsilon^2$
and $V(K)V(K^z)\leq V(C)V(C^o)$.
In particular, $C^o$ is an $o$-symmetric convex body with
axial rotational symmetry and satisfies
$\delta_{\rm BM}(C^o,B^n)>1+\gamma_1\varepsilon^2$.
According to Corollary~\ref{Falsestabcor},
there exist $u\in S^{n-1}$ and $a\in (0,\frac23)$
such that
$$
q(C^o(u,h_{C^o}(u)t))\geq 1+\gamma_2\varepsilon^6|\log\varepsilon|^{-1}
\mbox{ \ for $t\in(a,a+\gamma_2\varepsilon^6|\log\varepsilon|^{-1})$}.
$$
We may assume that $h_{C^o}(u)=1$.

Let $\widetilde{C}=C_{u^\bot}$.
Since the convexity of $C^o$ yields
 $|C^o(u,t)|\geq 4^{-(n-1)}|C^o(u,0)|$ if $t\leq \frac34$,
we deduce from  (\ref{Brunn}), (\ref{section}) and  (\ref{Brunnstab}) that
$$
V(\widetilde{C}^o)\geq 2\int_0^1|\widetilde{C}^o(u,t)|\,dt
\geq 2\int_0^1|C^o(u,t)|\,dt+
\gamma_3|C^o(u,0)|\varepsilon^{6n}|\log\varepsilon|^{-n}.
$$
On the one hand, $V(C^o)\leq 2|C^o(u,0)|$ by the Fubini Theorem 
and the Brunn-Minkowski inequality (\ref{Brunn}). Therefore,
\begin{eqnarray}
\nonumber
V(K)V(K^z)&\leq& V(C)V(C^o)\leq (1-\gamma_4\varepsilon^{6n}|
\log\varepsilon|^{-n})V(\widetilde{C})V(\widetilde{C}^o)\\
\label{BSstabeq}
&\leq &(1-\gamma_4\varepsilon^{6n}|\log\varepsilon|^{-n})\kappa_n^2,
\end{eqnarray}
which concludes the proof of Theorem~\ref{BSstab}.\\

\noindent{\bf Acknowledgement: } I am grateful for
the help of Keith M. Ball,  Andrea Colesanti, Matthieu Fradelizi,
Peter M. Gruber, Lars Hoffmann,
Greg Kuperberg, Erwin Lutwak, Endre Makai, Vitali D. Milman,
Shlomo Reisner, Paolo Salani, Rolf Schneider
in the preparation of this manuscript. I am also indebted
to the referee, whose remarks
substantially improved the presentation
of the paper.

\end{document}